\documentclass[12pt]{amsart}
\usepackage[left=1.5in, right=1.5in, top=1.5in, bottom=1.5in]{geometry}                % See geometry.pdf to learn the layout options. There are lots.
\usepackage[usenames,dvipsnames,svgnames,table]{xcolor}
\usepackage[pagebackref,linktocpage=true,colorlinks=true,linkcolor=Blue,citecolor=BrickRed,urlcolor=RoyalBlue]{hyperref}
\usepackage[alphabetic,msc-links,abbrev]{amsrefs} %shortalphabetic,msc-links,backrefs
\usepackage[colorinlistoftodos,prependcaption,textsize=tiny, textwidth=0.75in]{todonotes}

\newtheorem{thm}{Theorem}[section]

\theoremstyle{definition}

\newtheorem{rem}[thm]{Remark}
\usepackage{graphicx}
\usepackage{amssymb}
\usepackage{amsmath}
\usepackage{epstopdf}
\usepackage{hyperref}
\allowdisplaybreaks

\title[The equilibrium shape]{On the three-dimensional shape of a crystal}

\author{Emanuel Indrei and Aram Karakhanyan}
\address{Department of Mathematics and Statistics\\
Sam Houston State University\\
Huntsville, TX \\
USA.}
\address{School of Mathematics\\ The University of Edinburgh\\ Peter
Tait Guthrie
Road\\ EH9 3FD Edinburgh, UK.
}

\thanks{The research of the second author was partially supported
by EPSRC grant
EP/S03157X/1 {\em Mean curvature measure of free boundary}.}

\begin{document}
\setcounter{page}{1}
\pagenumbering{arabic}
\maketitle

\begin{abstract}
In this paper we completely settle the Almgren problem in $\mathbb R^3$ under some generic
conditions on the potential and tension functions.
The problem, among other things, appears in  classical thermodynamics when one is to understand if minimizing the free energy with convex potential and under a mass constraint generates a convex crystal.  
Our new idea in proving a three-dimensional convexity theorem is to utilize a stability theorem when $m$ is small, convexity when $m$ is small, and the first variation PDE with a new maximum principle approach.
\end{abstract}

\section{Introduction}
A fundamental problem in thermodynamics is to prove convexity of minimizers to the free energy minimization with mass constraint.
The free energy $\mathcal{E}(E)$ of a set of finite perimeter $E \subset \mathbb{R}^n$ with reduced boundary $\partial^* E$ is defined via the surface energy

$$
\mathcal{F}(E)=\int_{\partial^* E} f(\nu_E) d\mathcal{H}^{n-1},
$$
and, the potential energy 
$$
\mathcal{G}(E)=\int_E g(x)dx,
$$
where $g \ge 0$, $g(0)=0$:
$$
\mathcal{E}(E)=\mathcal{F}(E)+\mathcal{G}(E).     
$$
The following problem historically is attributed to Almgren.\\

{\bf Problem:} If the potential $g$ is convex (or, more generally, if the sub-level sets $\{g < t\}$
are convex), are minimizers convex or, at least, connected? \cite[p. 146]{MR2807136}.\\

In a recent paper, Indrei proved the 
existence of a convex $g \ge 0$, $g(0)=0$, so that there are no minimizers for $m>0$. Observe the general partition of the convexity problem into coercive (e.g. the monotone radial potential) and non-coercive potentials (e.g. the gravitational potential). Supposing $n=2$, under additional assumptions, the first author proved convexity for all $m>0$ \cite{Cryst} (cf. \cite{D}). 
In the argument, the planar context is crucial. Recently, the authors proved a sharp quantitative inequality for the isotropic radial Almgren problem ($f(x)=|x|, g=g(|x|)$) in $\mathbb{R}^n$. The theorem is the first positive result for all $m>0$ on the stability and convexity for a large class of potentials in higher dimension. 

For $g=0$ the stability appeared in 
\cite{MR2672283} with an explicit modulus; in
\cite{MR2456887} for
$g=0$ and the isotropic case with a semi-explicit modulus; and, in
\cite{Cryst} for
$m$ small with a semi-explicit modulus and a locally bounded
potential. 

Naturally, in physics, the most important dimension is $n=3$. 
We introduce a new method to prove:    
\begin{thm} \label{85}
If $g \in C^{2,\alpha_*}$ is convex, $f \in C^{4}(\mathbb{R}^3\setminus\{0\})$ is $\lambda-$elliptic, and 
$f,g$ admit minimizers $E_m \subset B_{R(m)}$ with $R \in L_{loc}^\infty(\mathbb{R}^+)$
either:\\
(i) $E_m$ is convex $\&$ unique for all $m \in (0,\infty)$;\\
(ii) there exist $\mathcal{M}>0$ $\&$ a modulus $w_m(0^+)=0$ such that for all $m \in (0, \mathcal{M})$, $E_m$ is unique, convex and there exist $\epsilon_0, \gamma>0$ such that for all $\epsilon \le \epsilon_0$,
$$
\liminf_{m \rightarrow \mathcal{M}^-} \frac{\gamma(\mathcal{M}-m)}{w_m(\epsilon)} \ge 1.
$$
\end{thm}
Our theorem implies convexity for a large collection of potentials; our argument is inclusive of also non-convex potentials. The main element is estimating the modulus.
\begin{rem}
If $f(\nu)=1$ when $\nu \in \mathbb{S}^{n-1}$, $g(x)=h(|x|)$, $h: \mathbb{R}^+ \rightarrow \mathbb{R}^+$ is increasing, $h(0)=0$,  then 
for any $m_1>m_0>0$, $\epsilon>0$,
$$
\inf_{m_1 \ge m \ge m_0} w_m(\epsilon)>0
$$
\cite{qk} and therefore 
for all $\mathcal{M}, \epsilon, \gamma>0$,
$$
\liminf_{m \rightarrow \mathcal{M}^-} \frac{\gamma(\mathcal{M}-m)}{w_m(\epsilon)} =0
$$
which precludes (ii). The result then yields uniqueness and convexity for any $m \in (0,\infty)$. 
\end{rem}
Our new idea is to utilize a stability theorem when $m$ is small, convexity when $m$ is small, and the first variation PDE with a new maximum principle approach. Assuming $g$ is coercive, the assumption on existence is true. Nevertheless, in certain configurations, one may prove existence for non-coercive potentials, e.g. the gravitational potential. \\
The stability result contains an invariance collection.
Define
\begin{align*}
\mathcal{A}_m=\mathcal{A}_{f,g,m}&=\{A: Ax=A_ax+x_a, x_a \in \mathbb{R}^3, \mathcal{E}(A_aE)=\mathcal{E}(E), \\
&|A_aE|=|E|=m \text{ for some minimizer $E$}  \}.
\end{align*}
An invariance map of the free energy is a transformation $A \in \mathcal{A}_m$.
The uniqueness of minimizers can only be true mod  $\mathcal{H}^{2}$ sets of measure zero and an invariance map generated by the mass, potential, and tension. In many classes of potentials, assuming $m$ is small, $A \in \mathcal{A}_m$ is a translation $Ax=x+z$, $z \in \mathbb{R}^3$. For example, suppose $g$ is zero on a ball $B$. If $m$ is small, note that uniqueness can only be shown up to a translation: $A_a=I_{3 \times 3}$, $x_a \in \mathbb{R}^3$ is such that $K_m+x_a \subset \{g=0\}$ when $K_m \subset B$ ($K_m$ is the Wulff shape such that $|K_m|=m$).  
The three transformations, reflection, rotation, and translation, always satisfy closure under convexity: $AE$ is convex iff $E$ is convex.

\section{Proof of Theorem \ref{85}} 
Define 
$$
\mathcal{A}_a=\{m: E_{\overline{m}} \hskip .08in \text{is unique $\&$ convex for all $0<\overline{m}\le m$}\}
$$
$$
\mathcal{M}=\sup \mathcal{A}_a.
$$
Theorem \ref{@'} and Theorem 2 in Figalli and Maggi \cite{MR2807136} imply $(0,m_a) \subset \mathcal{A}_a$. Hence $\mathcal{M} >0$. In addition, one may assume the invariance maps are closed under convexity. If $\mathcal{M}<\infty$, for $m \in (0,\mathcal{M})$, $E_m$ is unique $\&$ convex. 
Therefore either: (a) there exists a non-convex minimizer having mass $\mathcal{M}$; (b) there exist two convex minimizers not mod an invariance map equal having mass $\mathcal{M}$; or (c) for all $m \in (0, \mathcal{M}]$, $E_m$ is unique, convex and for $m>\mathcal{M}$ there exists $a<m$ such that either convexity or uniqueness fails for minimizers with mass $a$.
If $m_k<\mathcal{M}$, $m_k \rightarrow \mathcal{M}$, along a subsequence, $E_{m_k} \rightarrow T_\mathcal{M}$, with $|T_\mathcal{M}|=\mathcal{M}$, $T_\mathcal{M}$ a convex minimizer. 
Set 
$$
\epsilon=\frac{1}{5} \inf_{R} \frac{|R E_\mathcal{M} \Delta T_\mathcal{M}|}{|E_\mathcal{M}|}>0,
$$
where if (a) is valid, $E_\mathcal{M}$ is the non-convex minimizer and if (b) is true, $E_\mathcal{M}$ is a convex minimizer not (mod invariance transformations) equal to $T_\mathcal{M}$.
If $m \in (0,\mathcal{M})$, the uniqueness of convex minimizers  implies that there exists $w_m(\epsilon)>0$ such that for all $\epsilon>0$, if 
$|E|=|E_m|$, $E \subset B_R$, and 
$$
|\mathcal{E}(E)-\mathcal{E}(E_m)|<w_m(\epsilon),
$$
then there exists $R$ such that 
$$
\frac{|E_m \Delta RE|}{|E_m|}<\epsilon.
$$
Let $\{m_k\}$ be the sequence such that 
$$
\liminf_{m \rightarrow \mathcal{M}^-} \frac{\mathcal{M}^{\frac{2}{3}}-m^{\frac{2}{3}}}{w_m(\epsilon)}=\lim_{k \rightarrow \infty} \frac{\mathcal{M}^{\frac{2}{3}}-m_k^{\frac{2}{3}}}{w_{m_k}(\epsilon)},
$$
and define $\gamma_k$ via $|\gamma_k E_\mathcal{M}|=|E_{m_k}|$, i.e. $\gamma_k=(\frac{m_k}{\mathcal{M}})^{\frac{1}{3}}$. Note  

\begin{align*}
|\mathcal{E}(\gamma_k E_\mathcal{M})-\mathcal{E}(E_{m_k})|& \le |\mathcal{E}(\gamma_k E_\mathcal{M})-\mathcal{E}(E_\mathcal{M})|+|\mathcal{E}(T_\mathcal{M})-\mathcal{E}(E_{m_k})|\\
&\le\mathcal{F}(E_\mathcal{M})(1-\gamma_k^{2})+ (\sup_{B_R}g) |E_\mathcal{M} \Delta (\gamma_k E_\mathcal{M})|\\
&+|\mathcal{E}(T_\mathcal{M})-\mathcal{E}(E_{m_k})|.\\
\end{align*}
Moreover, 
\begin{align*}
\mathcal{E}(T_{\mathcal{M}})& \le \mathcal{E}(\frac{1}{\gamma_k}E_{m_k})\\
&=\frac{1}{\gamma_k^{2}}\mathcal{F}(E_{m_k})+ \int_{\frac{1}{\gamma_k}E_{m_k}}g(x)dx\\
&\le (\frac{1}{\gamma_k^{2}}-1)\mathcal{F}(E_{m_k})+(\sup_{B_R}g) |\frac{1}{\gamma_k}E_{m_k} \Delta E_{m_k}|+\mathcal{E}(E_{m_k})
\end{align*}
and similarly thanks to $|\frac{1}{\gamma_k}E_{m_k} \Delta E_{m_k}| \le a(\frac{1}{\gamma_k}-1)$ (e.g via \cite[Lemma 4]{MR2807136}) this implies  

$$
|\mathcal{E}(T_{\mathcal{M}})-\mathcal{E}(E_{m_k})| \le \alpha_p (\frac{1}{\gamma_k^{2}}-1)=\alpha (\mathcal{M}^{\frac{2}{3}}-m_k^{\frac{2}{3}}),
$$
$m_k \thickapprox \mathcal{M}$.

In particular, 
$$
|\mathcal{E}(\gamma_k E_\mathcal{M})-\mathcal{E}(E_{m_k})| \le \gamma_1 (\mathcal{M}^{\frac{2}{3}}-m_k^{\frac{2}{3}})
$$
where $\gamma_1=\gamma_1(\mathcal{M})$.

Suppose 

\begin{equation} \label{jw}
\liminf_{m \rightarrow \mathcal{M}^-} \frac{\mathcal{M}^{\frac{2}{3}}-m^{\frac{2}{3}}}{w_m(\epsilon)}<\frac{1}{\gamma_1}, 
\end{equation}
then for $k$ large
$$
|\mathcal{E}(\gamma_k E_\mathcal{M})-\mathcal{E}(E_{m_k})| \le  \frac{ \gamma_1 (\mathcal{M}^{\frac{2}{3}}-m_k^{\frac{2}{3}})}{w_{m_k}(\epsilon)} w_{m_k}(\epsilon)<w_{m_k}(\epsilon)
$$
and this implies the existence of $R_k$ such that 

$$
\frac{|E_{m_k} \Delta R_k(\gamma_k E_\mathcal{M})|}{|E_{m_k}|}<\epsilon.
$$
However, if $k$ is large, $\gamma_k \thickapprox 1$, which implies 
\begin{align*}
\frac{|(E_{m_k}) \Delta R_k(\gamma_k E_\mathcal{M})|}{|E_{m_k}|}& \thickapprox \frac{|T_\mathcal{M} \Delta R_k(E_\mathcal{M})|}{|E_\mathcal{M}|}\\
 &\ge \inf_{R} \frac{|RE_\mathcal{M} \Delta T_\mathcal{M}|}{|E_\mathcal{M}|} =5\epsilon,
\end{align*}
a contradiction.
Therefore \eqref{jw} is not true and 
$$
\liminf_{m \rightarrow \mathcal{M}^-} \frac{\mathcal{M}^{\frac{2}{3}}-m_k^{\frac{2}{3}}}{w_m(\epsilon)}\ge \frac{1}{\gamma_1}, 
$$
for $$\epsilon \le \epsilon_0:=\frac{1}{5}\inf_R \frac{|RE_\mathcal{M} \Delta T_\mathcal{M}|}{|E_\mathcal{M}|}.$$
Thus this yields $\gamma=\gamma(\mathcal{M})>0$, 
$$
\liminf_{m \rightarrow \mathcal{M}^-} \frac{\mathcal{M}-m_k}{w_m(\epsilon)}\ge \frac{1}{\gamma};
$$
observe the bound in (ii) is proved. The last part is to preclude (c). \\

\vskip .3in

Claim 1:
A convex minimizer at mass $\mathcal{M}$ is uniformly convex.
\vskip .3in

\noindent Proof of Claim 1:\\

The anisotropic mean curvature is 
$$
H_f=\text{trace} \Big(D^2 f A \Big),
$$
where $D^2f$ is the matrix of second tangential derivatives and $A$ is the second fundamental form. 
The formula for the first variation implies
\begin{equation} \label{fa}
H_f=\mu-g,
\end{equation}
where 
$$
\mu= \frac{2\mathcal{F}(E_{\mathcal{M}}) + \int_{\partial^* E_{\mathcal{M}}} g \langle x, \nu_{E_{\mathcal{M}}} \rangle d \mathcal{H}^{2}}{n|E_{\mathcal{M}}|}.
$$

Convexity of $E_{\mathcal{M}}$ and \eqref{fa} imply that locally there is a convex function $u\in C^{2.1}(\Omega), \Omega\subset\mathbb R^2$ so that  
$$a_{ij}(\nabla u)u_{ij}=\mu-g(x, u),$$
where 
$
a_{ij}(\nabla u)
,i,j \in \{1,2\}$, is a uniformly elliptic matrix given in terms of the 
second order derivatives of $f$ and depending on $\nabla u$
with $g$ being a convex function of $(x, u)\in \mathbb R^3$ and $\nabla$ the gradient, see Chapter 16.4 \cite{GT01}. 
Recall that for the classical case $f(\xi)=|\xi|$ we have 
$$
a_{ij}(\nabla u)=\frac1{\sqrt{1+|\nabla u|^2}}\left(\delta_{ij}-\frac{u_iu_j}{1+|\nabla u|^2}\right).
$$

Note
\begin{equation} \label{kx}
\mu-g>0 \hskip .1in \text{ on $\partial E_{\mathcal{M}}$}. 
\end{equation}
Indeed, let us choose a smoothly changing coordinate system 
so that $D^2u$ is diagonal. Then the mean curvature takes the form 
$H_f=a_{11}u_{11}+a_{22}u_{22}.$
After differentiating we get 
\begin{eqnarray*}
(H_f)_{ss}=- \nabla^2_{\mathbb R^3} g(x, u)
\partial_s\begin{pmatrix}
x_1\\x_2\\ u
\end{pmatrix}
\partial_s\begin{pmatrix}
x_1\\x_2\\ u
\end{pmatrix}
-g_u u_{ss}, \quad s=1, 2.
\end{eqnarray*}
Then 
\[
a_{ss}(H_f)_{ss}\le -g_u H_f,
\]
and consequently 
\[
a_{ss}(H_f)_{ss}-(g_u)^- H_f\le 0, 
\]
where $(g_u)^-$ is the negative part of $g_u$.
Hence the result follows from the strong minimum principle. 
\vskip .3in

\noindent Subclaim: If $\det D^2u(x_0)=0$ for some 
$x_0\in \Omega$ then $\det D^2u(x)=0$ for all  $x\in \Omega$.\\
\vskip .3in

\noindent Proof of Subclaim:\\
Observe that under our assumptions $u\in C^{3,1}(\Omega)$ thanks to  Corollary 16.7 \cite{GT01}.  
The proof is based on the observation that 
$w:=\det D^2u(x)$ satisfies an inequality of the form $a_{ij}w_{ij}-cw+b\cdot \nabla w\le 0$ near 
$x_0$, with $c\ge 0$.

Let us write the equation in the form 
$\mathcal F(D^2u, \nabla u)=\sum a_{ij}u_{ij}$ and 
let $f=\mu-g$, then the equation takes the form  
 
$$\mathcal F=f.$$
 Differentiate twice in $x_s, x_t, 1\le s, t\le2$ to get 
 \begin{eqnarray*}
\mathcal F_s=f_s,\\
 \mathcal F_{st}=f_{st}.
 \end{eqnarray*}

Now we have that 
\begin{eqnarray*}
w_s=u^{ij}u_{ijs}, \quad w_{st}=u^{ij, kl}u_{ijs}u_{klt}+u^{ij}u_{ijst},
\end{eqnarray*}
where $u^{ij}$ is the cofactor matrix.

On the other hand 
\begin{eqnarray*}
\mathcal F_s
&=&
\nabla _p a_{lm}\cdot \nabla u_{s} u_{lm}+a_{lm}u_{lms},\\
%%%
\mathcal F_{st}
&=&
(\nabla_{pp}^2 a_{lm}\nabla u_{t})\cdot \nabla u_{s}u_{lm}\\
&&
+
\nabla_p a_{lm}\cdot \nabla u_{st}u_{lm}+\nabla_pa_{lm}\cdot \nabla u_{s}u_{lmt}
+
\nabla_p a_{lm}\cdot \nabla u_{t}u_{lms}\\
&&
+
a_{lm}u_{lmst}\\
&:=&
\mathcal F^{(2)}+ \mathcal F^{(3)}+
a_{lm}u_{lmst},
\end{eqnarray*}
where we use the notation with dummy variable
$p:=\nabla u$.

Since the Weingarten mapping is self-adjoint, then at each point $x$, near $x_0$
we have 
\begin{equation}\label{eq:diag}
D^2 u(x)=\mbox{diag}[\lambda_1, \lambda_2]
\end{equation} 
in a continuously changing  coordinate system.
Moreover, $\lambda_2\ge\lambda_1\ge 0$. By \eqref{kx}, $\lambda_1+\lambda_2>0$.
Suppose $w(x_0)=0$, then $u_{11}(x_0)u_{22}(x_0)=0$ and without loss of generality 
\begin{equation}\label{eq:delta}
u_{11}(x)>10\delta\quad \mbox{and}\quad u_{22}(x)<\delta
\end{equation}
for $\delta>0$, 
in some neighborhood $x\in B_{r_0}(x_0)$, $r_0>0$ small.
Using these observations we can make the following explicit computations
\begin{eqnarray}
w_s&=&u_{11s}u_{22}+u_{11}u_{22s},\\
w_{st}&=&u_{11st}u_{22}+u_{11s}u_{22t}+u_{11t}u_{22s}+u_{11}u_{22st}-2u_{12t}u_{12s}.
\end{eqnarray}
The second order derivatives appearing in $\mathcal F_{st}$, after contracting 
with the cofactor matrix $u^{ij}=\mbox{diag}(u_{22}, u_{11})$,  and using \eqref{eq:diag}, can be simplified as follows
\begin{eqnarray*}
u^{st}\mathcal F_{st}^{(2)}
&:=&
(\nabla_{pp}^2 a_{lm}\nabla u_{t})\cdot \nabla u_{s}u_{lm}u^{st}\\
&=&
(\nabla_{pp}^2 a_{lm}\nabla u_{1})\cdot \nabla u_{1}u_{lm}u^{11}+
(\nabla_{pp}^2 a_{lm}\nabla u_{2})\cdot \nabla u_{2}u_{lm}u^{22}\\
&=&
(\nabla_{pp}^2 a_{lm}\nabla u_{1})\cdot \nabla u_{1}u_{lm}u_{22}+
(\nabla_{pp}^2 a_{lm}\nabla u_{2})\cdot \nabla u_{2}u_{lm}u_{11}\\
&=&
(\nabla_{pp}^2 a_{11}\nabla u_{1})\cdot \nabla u_{1}u_{11}u_{22}+
(\nabla_{pp}^2 a_{22}\nabla u_{1})\cdot \nabla u_{1}u_{22}u_{22}\\
&&
+
(\nabla_{pp}^2 a_{11}\nabla u_{2})\cdot \nabla u_{2}u_{11}u_{11}
+
(\nabla_{pp}^2 a_{22}\nabla u_{2})\cdot \nabla u_{2}u_{22}u_{11}\\
&=&
\left((\nabla_{pp}^2 a_{11}\nabla u_{1})\cdot \nabla u_{1}+ (\nabla_{pp}^2 a_{22}\nabla u_{2})\cdot \nabla u_{2}\right)w\\
&&
+
(\partial_{p_1p_1}a_{22}+\partial_{p_2p_2}a_{11})w^2.
\end{eqnarray*}
Consequently, 
\begin{equation}\label{eq:F2}
u^{st}\mathcal F_{st}^{(2)}=O(cw+b\cdot \nabla w), 
\end{equation}
for some fixed $c>0$ and $b\in \mathbb R^2$.

Next, let us compute the expression 
\begin{eqnarray}\nonumber
a_{lm}w_{lm}
&=&
u^{st}a_{lm}u_{imst}+a_{lm}u^{st, ij}u_{stl}u_{ijm}\\\nonumber
&=&
-u^{st}\mathcal F_{st}^{(2)}-u^{st}\mathcal F_{st}^{(3)}+u^{st}f_{st}
+a_{lm}(u_{11l}u_{22m}+u_{11m}u_{22l}-2u_{12m}u_{12l})\\\label{eq:linearized}
&=&
-u^{st}\mathcal F_{st}^{(2)}-u^{st}\mathcal F_{st}^{(3)}+u^{st}f_{st}
+J.
\end{eqnarray}
We need to simplify the last term 
$J:=a_{lm}(u_{11l}u_{22m}+u_{11m}u_{22l}-2u_{12m}u_{12l})$.
It can be written in a more explicit form as follows
\begin{eqnarray*}
J
&=&
a_{11}(u_{111}u_{221}+u_{111}u_{221}-2u_{121}^2)
+
a_{12}(u_{111}u_{222}+u_{112}u_{221}-2u_{122}u_{121})\\
&&
+
a_{21}(u_{112}u_{221}+u_{111}u_{222}-2u_{121}u_{122})
+
a_{22}(u_{112}u_{222}+u_{112}u_{222}-2u_{122}^2)\\
&=&
2\left(a_{11}(u_{111}u_{221}-u_{121}^2)
+
a_{12}(u_{111}u_{222}-u_{122}u_{121})
+
a_{22}(u_{112}u_{222}-u_{122}^2)
\right).
\end{eqnarray*}
Using the explicit forms of $w_s, \mathcal F_s$ 
we obtain
\begin{eqnarray}\label{eq:mek}
a_{11}u_{11s}+2a_{12}u_{12s}+a_{22}u_{22s}=f_s-(\partial_{p_s}a_{ll})u_{ss}u_{ll},\\
\label{eq:erku}
u_{11s}u_{22}+u_{11}u_{22s}=w_s,
\end{eqnarray}
since 
\[
\nabla_p a_{lm}\cdot \nabla u_su_{lm}
=(\partial_{p_s}a_{ll})u_{ss}u_{ll}.
\]
From \eqref{eq:erku} 
\begin{equation}\label{eq:fff}
u_{22s}=\frac{w_s-u_{11s}u_{22}}{u_{11}}, 
\end{equation} 
plugging this into \eqref{eq:erku}
yields
\begin{eqnarray*}
f_s-(\partial_{p_s}a_{ll})u_{ss}u_{ll}
&=&
a_{11}u_{11s}+2a_{12}u_{12s}+a_{22}\frac{w_s-u_{11s}u_{22}}{u_{11}}\\
&=&
u_{11s}\frac{a_{11}u_{11}-a_{22}u_{22}}{u_{11}} + a_{22}\frac{w_s}{u_{11}}
+
2a_{12}u_{12s}.
\end{eqnarray*}
If $s=1$, then $u_{12s}=u_{121}=u_{112}$, and from the above computation 
\begin{eqnarray*}
f_1-(\partial_{p_1}a_{ll})u_{11}u_{ll}=u_{111}\frac{a_{11}u_{11}-a_{22}u_{22}}{u_{11}} + a_{22}\frac{w_1}{u_{11}}
+
2a_{12}u_{112}.
\end{eqnarray*}
Similarly, for $s=2$ we obtain 
\begin{eqnarray*}
f_2-(\partial_{p_2}a_{ll})u_{22}u_{ll}
=
u_{112}\frac{a_{11}u_{11}-a_{22}u_{22}}{u_{11}} + a_{22}\frac{w_2}{u_{11}}
+
2a_{12}\frac{w_1-u_{111}u_{22}}{u_{11}}.
\end{eqnarray*}
Combining the last two equation we get a system of equations for the remaining third order derivatives $u_{111}$
and $u_{112}$;
\begin{eqnarray*}
f_1-(\partial_{p_1}a_{ll})u_{11}u_{ll}-a_{22}\frac{w_1}{u_{11}}=u_{111}\frac{a_{11}u_{11}-a_{22}u_{22}}{u_{11}}  
+
2a_{12}u_{112},\\
f_2-(\partial_{p_2}a_{ll})u_{22}u_{ll} -a_{22}\frac{w_2}{u_{11}} -2a_{12}\frac{w_1}{u_{11}}
=
-
2a_{12}\frac{u_{22}}{u_{11}} u_{111}+u_{112}\frac{a_{11}u_{11}-a_{22}u_{22}}{u_{11}} .
\end{eqnarray*}
Note that the determinant of the coefficient matrix is 
\[
\mathcal D:=\frac{(a_{11}u_{11}-a_{22}u_{22})^2}{u_{11}^2}+4a_{12}^2\frac{u_{22}}{u_{11}} >0, 
\]
and, moreover,
\begin{equation}\label{eq:discriminant}
\frac1{\mathcal D}= \frac{u_{11}^2}{(a_{11}u_{11}-a_{22}u_{22})^2}\left(1+4a_{12}^2\frac{u_{22}}{u_{11}}h\right)
\end{equation}
for some bounded  function $h$ in view of \eqref{eq:delta}.

Solving the system we find
\begin{eqnarray*}
u_{111}
&=&
\frac1{\mathcal D}
\Big(\frac{a_{11}u_{11}-a_{22}u_{22}}{u_{11}}(f_1-(\partial_{p_1}a_{ll})u_{11}u_{ll}-a_{22}\frac{w_1}{u_{11}})\\
&&-
2a_{12}(f_2-(\partial_{p_2}a_{ll})u_{22}u_{ll} -a_{22}\frac{w_2}{u_{11}} -2a_{12}\frac{w_1}{u_{11}})
\Big)\\
&=&
\frac1{\mathcal D}
\Big(\frac{a_{11}u_{11}-a_{22}u_{22}}{u_{11}}(f_1-(\partial_{p_1}a_{11})u_{11}^2)
-
2a_{12}(f_2-(\partial_{p_2}a_{22})u_{22}^2)\Big)
+
O(cw+b\cdot \nabla w)\\
&=&
\frac1{\mathcal D}
\Big(\frac{a_{11}u_{11}-a_{22}u_{22}}{u_{11}}(f_1-(\partial_{p_1}a_{11})u_{11}^2)
-
2a_{12}(f_2)\Big)
+
O(cw+b\cdot \nabla w)
\end{eqnarray*}
and
\begin{eqnarray*}
u_{112}
&=&
\frac1{\mathcal D}
\Big(
\frac{a_{11}u_{11}-a_{22}u_{22}}{u_{11}}(f_2-(\partial_{p_2}a_{ll})u_{22}u_{ll} -a_{22}\frac{w_2}{u_{11}} -2a_{12}\frac{w_1}{u_{11}})\\
&&+
2a_{12}\frac{u_{22}}{u_{11}}(f_1-(\partial_{p_1}a_{ll})u_{11}u_{ll}-a_{22}\frac{w_1}{u_{11}})
\Big)\\
% last line
&=&
\frac1{\mathcal D}
\Big(
\frac{a_{11}u_{11}-a_{22}u_{22}}{u_{11}}(f_2-(\partial_{p_2}a_{22})u_{22}^2)+
2a_{12}\frac{u_{22}}{u_{11}}(f_1-(\partial_{p_1}a_{11})u_{11}^2)
\Big)+
O(cw+b\cdot \nabla w)\\
&=&
\frac1{\mathcal D}
\Big(
\frac{a_{11}u_{11}-a_{22}u_{22}}{u_{11}}f_2
\Big)+
O(cw+b\cdot \nabla w)\\
&=&
\frac{u_{11}}{a_{11}u_{11}-a_{22}u_{22}}\left(1+4a_{12}^2\frac{u_{22}}{u_{11}}h\right)f_2
+
O(cw+b\cdot \nabla w)\\
&=&
\frac{u_{11}}{a_{11}u_{11}-a_{22}u_{22}}f_2
+
O(cw+b\cdot \nabla w)
.
\end{eqnarray*}
Therefore, combining with \eqref{eq:discriminant} and \eqref{eq:delta} we infer  that 
$u_{112}$ and $u_{112}$ can be estimated in terms of the lower order derivatives of $u$, hence we conclude that 
\begin{equation}\label{eq:system-solve}
u_{111}, u_{112}=O(cw+b\cdot \nabla w).
\end{equation}

Returning to 
\begin{eqnarray*}
J
&=&
2\Bigg\{a_{11}\left(u_{111}\frac{w_1-u_{111}u_{22}}{u_{11}}-u_{112}^2\right)\\
&&
+
a_{12}\left(u_{111}\frac{w_2-u_{112}u_{22}}{u_{11}}-\frac{w_1-u_{111}u_{22}}{u_{11}}u_{112}\right)\\
&&
+
a_{22}\left(u_{112}\frac{w_2-u_{112}u_{22}}{u_{11}}-\left(\frac{w_1-u_{111}u_{22}}{u_{11}}\right)^2\right)
\Bigg\}\\
&=&
2\Bigg\{
u_{111}^2(-a_{11}\frac{u_{22}}{u_{11}}-a_{22}\frac{u_{22}^2}{u_{11}^2})
+
u_{112}^2(-a_{11}-a_{22} \frac{u_{22}}{u_{11}})\\
&&
+
u_{111}(w_1 \frac{a_{11}}{u_{11}}+ w_2 \frac{a_{22}}{u_{11}}+2w_1 \frac{a_{22}u_{22}}{u_{11}})\\
&&
+
u_{112}(-w_1 \frac{a_{12}}{u_{11}}+w_2 \frac{a_{22}}{u_{11}})-a_{22}\frac{w_1^2}{u_{11}^2}
\Bigg\}\\
&\le&
u_{111}(w_1 \frac{a_{11}}{u_{11}}+ w_2 \frac{a_{22}}{u_{11}}+2w_1 \frac{a_{22}u_{22}}{u_{11}})\\
&&
+
u_{112}(-w_1 \frac{a_{12}}{u_{11}}+w_2 \frac{a_{22}}{u_{11}})\\
&=&
O(cw+b\cdot \nabla w), 
\end{eqnarray*}
where the last line follows from \eqref{eq:system-solve} and \eqref{eq:delta}.

For the third order derivatives in $\mathcal F_{st},$ after contraction with 
$u^{ij}$ we have
\begin{eqnarray*}
u^{st}\mathcal F_{st}^{(3)}
&:=&
\left(\nabla_p a_{lm}\cdot \nabla u_{st}u_{lm}+\nabla_pa_{lm}\cdot \nabla u_{s}u_{lmt}
+
\nabla_p a_{lm}\cdot \nabla u_{t}u_{lms}\right)u^{st}\\
&=&
\sum_{s=1}^2 \sum_{l,m}\left( \nabla_p a_{lm}\cdot \nabla u_{ss}u_{lm}u^{ss}
+
2\nabla_pa_{lm}\cdot \nabla u_{s}u_{lms}u^{ss}\right)\\
&=&
\sum_{s=1}^2\left( \nabla_p a_{11}\cdot \nabla u_{ss}u_{11}u^{ss}
+
\nabla_p a_{22}\cdot \nabla u_{ss}u_{22}u^{ss}
+
2\sum_{l,m}\nabla_pa_{lm}\cdot \nabla u_{s}u_{lms}u^{ss}\right)\\
%%%
&=&
\sum_{s=1}^2\left( \partial_{p_1} a_{11} u_{1ss}u_{11}u^{ss}
+\partial_{p_2} a_{11}u_{2ss}u_{11}u^{ss}\right)\\
&&
+
\sum_{s=1}^2\left( \partial_{p_1} a_{22}  u_{1ss}u_{22}u^{ss}
+  \partial_{p_2} a_{22}  u_{2ss}u_{22}u^{ss}\right)\\
&&
+
2\sum_{s=1}^2\sum_{l,m}\left( \partial_{p_1}a_{lm}u_{1s}u_{lms}u^{ss}
+\partial_{p_2}a_{lm}u_{2s}u_{lms}u^{ss}\right)\\
&=&
(\partial_{p_1} a_{11}u_{11}) w_1+(\partial_{p_2} a_{11}u_{11})w_2\\
&&
+
(\partial_{p_1} a_{22} u_{22})w_1+ (\partial_{p_2} a_{22} u_{22})w_2\\
&&
+
2\sum_{l,m}\left( \partial_{p_1}a_{lm}u_{lm1}
+\partial_{p_2}a_{lm}u_{lm2}\right)w.
\end{eqnarray*}

From here and our estimates for the third order derivatives we conclude that 
\begin{equation}\label{eq:F3}
u^{st}\mathcal F_{st}^{(3)}=O(cw+b\cdot \nabla w), 
\end{equation}
for some fixed $c>0$ and $b\in \mathbb R^2$.

Using this and \eqref{eq:linearized} we get 
\begin{eqnarray}\label{eq:rox}
a_{lm}w_{lm}\le u^{st}f_{st}+ O(cw+b\cdot \nabla w).
\end{eqnarray}
To finish the proof note that 
\begin{eqnarray*}
u^{st}f_{st}
&=&
u^{st}\nabla_{\mathbb R^3}f(x, u)\cdot \begin{pmatrix}
0\\0\\u_{st}
\end{pmatrix}
+
u^{st}\nabla^2_{\mathbb R^3} f(x, u)
\partial_s\begin{pmatrix}
x_1\\x_2\\ u
\end{pmatrix}
\partial_t\begin{pmatrix}
x_1\\x_2\\ u
\end{pmatrix}\\
&=&
f_u w
-
u^{ss}\nabla^2_{\mathbb R^3} g(x, u)
\partial_s\begin{pmatrix}
x_1\\x_2\\ u
\end{pmatrix}
\partial_s\begin{pmatrix}
x_1\\x_2\\ u
\end{pmatrix}\\
&\le&
f_u w
\end{eqnarray*}
since we assume that $g$ is convex.
Summarizing, it follows from the last inequality and \eqref{eq:rox} that
\[
a_{lm}w_{lm} + c w+b\cdot \nabla w\le 0.
\] 
Writing $c=c^+-c^-, c^\pm\ge 0$, and using $w\ge 0$ we get that 
\begin{equation}
a_{lm}w_{lm}- c^- w+b\cdot \nabla w\le 0.
\end{equation}
Applying the strong minimum principle we see that 
$w=0$ in $B_{r_0}(x_0)$. Therefore, the proof of Subclaim is finished.

Next, we prove Claim 1: if the Gauss curvature of $\partial E_{\mathcal{M}}$ vanishes at some point, $w(x_0)=0$, then the Gauss curvature is zero everywhere on $\partial E_{\mathcal{M}}$ (Subclaim).
By Theorem 2.8 \cite{RT} $u$ is the lower boundary of the convex hull of the set of points
$(x, u|_{\partial\Omega}),$ for any strictly convex $\Omega$. 
For such $\Omega,$ if we pick a point $x\in \Omega$ then there is a line segment passing 
through $x$. These line segments cannot intersect since otherwise that mean curvature vanishes at the 
intersection. Thus the graph of $u$ over $\Omega$ is a ruled surface. 
If we take a hyperplane perpendicular to the one containing the domain $\Omega$, then 
for $\Omega^\perp$ lying on this hyperplane  
the same conclusion will hold. However,
 the line segments 
generated by $\Omega$ and $\Omega^\perp$ must intersect, which will contradict the $C^3$ regularity of the surface. This yields the proof of Claim 1.\\

Suppose for $m_k>\mathcal{M}$ there is $m_{j_k}<m_k$ with $E_{m_{j_k}}$ a non-convex minimizer. Via Claim 1, $E_\mathcal{M}$ is uniformly convex. In particular, the two curvatures are uniformly positive. Via the smoothness, up to a subsequence,  $E_{m_{j_k}} \rightarrow E_{\mathcal{M}}$ in $C^2$. Observe that for $k$ sufficiently large, the regularity implies that the principal curvatures of  $E_{m_{j_k}}$ are near the ones of $E_{\mathcal{M}}$ and thus this contradicts non-convexity. In particular, 
\begin{equation} \label{y*}
\text{if $m>\mathcal{M}$ is near $\mathcal{M}$, then $E_m$ is uniformly convex.} 
\end{equation}
To show uniqueness the next fact is sufficient:\\ 

\noindent The Uniqueness Fact: There exists $m_0>0$ and a modulus of continuity $a(m,0^+)=0$ such that for all $m<\mathcal{M}+m_0$ there exists $\epsilon_0>0$ such that for all $0<\epsilon<\epsilon_0$ $\&$ for all minimizers $E_m \subset B_R$, $E \subset B_R$, $|E|=|E_m|=m<\mathcal{M}+m_0$, if 
$$
|\mathcal{E}(E_m)-\mathcal{E}(E)| < a(m,\epsilon),
$$
there exists an invariance map $A$  such that  
$$
\frac{|E \Delta A E_m|}{|E_m|} < \epsilon.
$$
Assume the uniqueness is false. Then for all $m_0>0$, for all moduli $q$ there exists $m<\mathcal{M}+m_0$ such that for a fixed $\epsilon_0 >0$ there exists $\epsilon<\epsilon_0$ $\&$ there exist $E_{m, \epsilon_0}, E_{m, \epsilon_0}' \subset B_R$, $|E_{m, \epsilon_0}|=|E_{m, \epsilon_0}'|=m$ such that

$$
|\mathcal{E}(E_{m,\epsilon_0})-\mathcal{E}(E_{m,\epsilon_0}')| < q_m(\epsilon),
$$ 
and

\begin{equation} \label{wl4}
\inf_{A} \frac{|E_{m,\epsilon_0}' \Delta AE_{m,\epsilon_0}|}{|E_{m,\epsilon_0}|}  \ge \epsilon>0.
\end{equation}
Let $m_0=\frac{1}{k}$, $w_k \rightarrow 0^+$, $\hat{q}$ a modulus of continuity and define 

\begin{equation} \label{t3}
q_k=w_k \hat{q}(\epsilon),
\end{equation}
hence there exists $m_k<\mathcal{M}+\frac{1}{k}$ such that for a fixed $\epsilon_0 >0$ there exists $\epsilon<\epsilon_0$ $\&$ there exist minimizers $E_{m_k, \epsilon_0}$, in addition some sets $E_{m_k, \epsilon_0}' \subset B_R$, $|E_{m_k, \epsilon_0}|=|E_{m_k, \epsilon_0}'|=m_k<\mathcal{M}+\frac{1}{k}$ such that
$$
|\mathcal{E}(E_{m_k,\epsilon_0})-\mathcal{E}(E_{m_k,\epsilon_0}')| < q_k,
$$ 
and
\begin{equation} \label{wlje1}
\inf_{A} \frac{|E_{m_k,\epsilon_0}' \Delta AE_{m_k,\epsilon_0}|}{|E_{m_k,\epsilon_0}|} \ge \epsilon>0.
\end{equation}
Set

$E_{m_k}=E_{m_k,\epsilon_0}$, $E_{m_k}'=E_{m_k,\epsilon_0}'$. Also, define $\gamma_k=(\frac{\mathcal{M}}{m_k})^{\frac{1}{3}}$ such that
$$|\gamma_k E_{m_k}|=|E_\mathcal{M}|.$$  \\
Next, observe that thanks to the compactness, $E_{m_k, \epsilon_0} \rightarrow E$, this yields $\mathcal{E}(E_{m_k, \epsilon_0}) \rightarrow \mathcal{E}(E)$, where $E$ is a minimizer, $|E|=\mathcal{M}$. In addition,
$$
|\mathcal{E}(E_{m_k,\epsilon_0})-\mathcal{E}(E_{m_k,\epsilon_0}')| < q_k \rightarrow 0
$$ 
also implies along a subsequence
$$E_{m_k,\epsilon_0}' \rightarrow \hat{E},$$
$|\hat{E}|=\mathcal{M},$ $\hat{E}$ a minimizer. 
The aforementioned
\begin{equation} \label{wlje1}
\inf_{A} \frac{|E_{m_k,\epsilon_0}' \Delta AE_{m_k,\epsilon_0}|}{|E_{m_k,\epsilon_0}|} \ge \epsilon>0
\end{equation}
therefore yields a contradiction:  initially, the uniqueness at mass $\mathcal{M}$ yields $A_1$ so that $A_1 E=\hat{E}$; thus
$$
A_1E_{m_k,\epsilon_0} \rightarrow A_1 E,
$$
$$
E_{m_k,\epsilon_0}'  \rightarrow \hat{E},
$$

$$
\frac{|E_{m_k,\epsilon_0}' \Delta A_1E_{m_k,\epsilon_0}|}{|E_{m_k,\epsilon_0}|} \rightarrow 0.
$$
Hence \eqref{y*} together with uniqueness preclude (c).

\section{Appendix} \label{tl;j}

\subsection{Modulus of the free energy}
If $g$ is locally bounded, the subsequent theorem solves the more general uniqueness problem in any dimension.
\begin{thm}[\cite{Cryst}] \label{@'}
Suppose $g \in L_{loc}^\infty(\{g<\infty\})$ admits minimizers $E_m \subset B_{R}$ for all $m$ small. There exists $m_0>0$ and a modulus of continuity $q(0^+)=0$ such that for all $m<m_0$ there exists $\epsilon_0>0$ such that for all $0<\epsilon<\epsilon_0$ and for all minimizers $E_m \subset B_R$, $E \subset B_R$, $|E|=|E_m|=m<m_0$, if 
$$
|\mathcal{E}(E_m)-\mathcal{E}(E)| < a(m,\epsilon)=q(\epsilon)m^{\frac{n-1}{n}},
$$
there exists an invariance map $A \in \mathcal{A}_m$  such that  
$$
\frac{|E \Delta A E_m|}{|E_m|} < \epsilon.
$$
Also, $AE_m \approx E_m + \alpha_m$: there exists $\alpha_m \in \mathbb{R}^n$, $c(n)>0$, so that 
$$
|AE_m\Delta \Big(E_m + \alpha_m \Big)| \le 2 \Big(\frac{1}{c(n)} \frac{1}{n|K|^{\frac{1}{n}}} (\sup_{B_{R_m}} g)\Big)^{\frac{1}{2}} m^{1+\frac{1}{2n}},
$$
where the radius $R_m>0$ is such that 
$$
(\frac{m}{|K|})^{\frac{1}{n}}K \subset B_{R_m}.
$$
\end{thm}

%\newpage

\end{document}